\documentclass[12pt,a4paper]{article}

\usepackage{amsmath}
\usepackage{amsthm}
\usepackage{amssymb}
\usepackage{cain_arxiv}
\usepackage{hyperref}

\theoremstyle{definition}

\theoremstyle{plain}

\newtheorem*{maintheorem}{Main theorem}
\newtheorem*{lemma*}{Lemma}

\numberwithin{equation}{section}

\DeclareMathOperator{\gL}{\mathcal{L}}
\DeclareMathOperator{\gR}{\mathcal{R}}
\DeclareMathOperator{\gD}{\mathcal{D}}
\DeclareMathOperator{\gJ}{\mathcal{J}}

\newcommand{\imreduces}{\rightarrow}

\DeclareMathOperator{\relrho}{\mathnormal{\rho}}

\begin{document}

\title{A simple non-bisimple congruence-free finitely presented monoid}
\author{Alan J. Cain \& Victor Maltcev}
\date{}

\thanks{The first author's research was funded by the European
  Regional Development Fund through the programme {\sc COMPETE} and by
  the Portuguese Government through the {\sc FCT} (Funda\c{c}\~{a}o
  para a Ci\^{e}ncia e a Tecnologia) under the project {\sc
    PEst-C}/{\sc MAT}/{\sc UI0}144/2011 and through an {\sc FCT}
  Ci\^{e}ncia 2008 fellowship.}

\maketitle

\address[AJC]{%
Centro de Matem\'{a}tica, Universidade do Porto, \\
Rua do Campo Alegre 687, 4169--007 Porto, Portugal
}
\email{%
ajcain@fc.up.pt
}
\webpage{%
www.fc.up.pt/pessoas/ajcain/
}

\address[VM]{%
Mathematical Institute, University of St Andrews,\\
North Haugh, St Andrews, Fife KY16 9SS, United Kingdom
}
\email{%
victor.maltcev@gmail.com
}

\begin{abstract}
We exhibit an example of a finitely presented monoid that is
congruence-free and simple but not bisimple.

\end{abstract}

\dummysection

There are several known examples of finitely presented infinite simple
groups, and of their semigroup theory counterparts, finitely presented
infinite congruence-free semigroups. Because of Rees congruences, it
is immediate that every congruence-free monoid is either simple or
$0$-simple. Within the class of finitely presented semigroups, there
are known examples of infinite congruence-free monoids that are
$0$-simple~\cite{acmu_anb,maltcev_finite}, and that are
simple~\cite{acmu_anb,maltcev_finite}. However, all previously known
examples of finitely presented congru\-ence-free \emph{simple} monoids
were also \emph{bisimple}. In an earlier paper, we asked whether there
existed a finitely presented congruence-free mon\-oid that was not
bisimple \cite[Open Problem 5.1]{acmu_anb}. The aim of this note is to
answer this question positively by exhibiting an example showing that
bisimplicity is not necessarily a consequence of simplicity,
congruence-freeness, and finite presentability. For further background
on finitely presented and more generally finitely generated
congruence-free semigroups we refer the reader
to~\cite{birget_monoid,birget_rl,birget_bernoulli,birget_monoidsthatmap,birget_thompson,byleen_embedbisimple,byleen_embednoidem,byleen_embedcongfree}.

\begin{maintheorem}
The monoid $M$ presented by the rewriting system
\begin{align*}
e^3 &\imreduces e,\\
xey &\imreduces y,\\
xe^2y &\imreduces x,\\
xy &\imreduces 1,
\end{align*}
over the alphabet $A = \{x,y,e\}$ is congruence-free and simple but not bisimple.
\end{maintheorem}

\setprooftag{Theorem}

\begin{proof}
First of all, note that the rewriting system is complete: there are no
overlaps of left-hand sides of rewriting rules, which implies the
system is confluent, and every rewriting rule descreases length, which
implies the system is noetherian.

The language of normal forms for the given rewriting system is 
\[
\{e,y\}^*\{e,x\}^* - A^*e^3A^*.
\]
It is easy to see that the normal forms of right-invertible elements
of $M$ are words from $\{1\}\cup x\{e,x\}^* - A^*e^3A^*$, and the
normal forms of left-invertible elements of $M$ are words from
$\{1\}\cup\{e,y\}^*y - A^*e^3A^*$. In particular it follows that for
every element $w\in M$ there exist $p,q\in M$ such that
$pwq\in\{1,e,e^2\}$. Since $xe^2y^2 = 1$ it follows that every element
is $\gJ$-related to $1$ and so $M$ is simple. But $M$ is not bisimple
since $(e,1)\notin\gD$. Indeed, $e\gD 1$ would imply
that $e$ is a product of a left-invertible element by a
right-invertible element, which is impossible.

It remains to prove that $M$ is congruence-free; that is, to prove
that the only congruences on $M$ are the equality relation $=$ and the
universal relation $M \times M$. We will proceed by complete
induction: for two distinct normal form words $u,v\in M$, we will
prove that if $\relrho$ is a congruence on $M$ such that $u\relrho v$,
then $\relrho=M\times M$; the induction will be taken on $|u|+|v|$.

The base case is $|u|+|v|=1$, i.e. when one of $u$ and $v$ is $1$ and
the other from $A$:
\begin{itemize}
\item If $e \relrho 1$, then $xey \relrho xy$ and so $y = xey \relrho
  xy = 1$ and hence $x \relrho xy = 1$; since $A$ generates $M$, we
  have $\relrho =M \times M$.
\item If $x \relrho 1$, then $y \relrho xy = 1$ and so $e \relrho xey
  = y \relrho 1$; again $\relrho = M \times M$.
\item If $y \relrho 1$, then $x \relrho xy = 1$ and so $e \relrho 1$
  as in the previous case and again $\relrho = M \times M$.
\end{itemize}

Now we prove the induction step. Before we begin, we prove the
following auxiliary result.

\begin{lemma*}
Every group homomorphic image of $M$ is trivial, and hence the only group
congruence on $M$ is $M \times $M.
\end{lemma*}

\setprooftag{Lemma}

\begin{proof}[Proof]
Let $N$ be a group homomorphic image of $M$. Then from $xey=y$ it
follows that $x=e^{-1}$ in $N$, and from $xe^2y=x$ it follows that
$y=e^{-2}$ in $N$. Then $e^{-1}=e^{-3}=xy=1$ and so $x=y=e=1$ in $N$.
\end{proof}

Returning to the induction step, let $u$ and $v$ be as in the
hypothesis and suppose the result holds for all pairs of normal form
words the sum of whose lengths is strictly less than $|u|+|v|$. We
consider five cases: both $u$ and $v$ contain letters $x$; only one
contains a letter $x$; both contain letters $y$ but not $x$; neither
contains letters $x$ and only one contains letters $y$; and finally
neither word contains letters $x$ or $y$.

\begin{enumerate}
\item Both $u$ and $v$ contain letters
$x$. Then $u\equiv Uxe^{\alpha}$ and $v\equiv Vxe^{\beta}$, where
$\alpha,\beta\in\{0,1,2\}$. If $\alpha=\beta$, then $U \neq V$ (since
$u \neq v$) and $U=ue^{2-\alpha}y^2\relrho ve^{2-\alpha}y^2=V$, and so
by induction $\rho = M \times M$ and the proof is complete. So we may
assume that $\alpha\neq\beta$. Interchanging $u$ and $v$ if necessary,
assume that $\alpha<\beta$. Then we have three sub-cases to consider:
\begin{enumerate}
\item $\alpha=1$ and $\beta=2$. Then $Uxe\relrho Vxe^2$. Hence
  $Uy=Uxey\relrho Vxe^2y=Vx$ and so $Uy^2\relrho Vxy=V$. Also we have
  $Uxe^2\relrho Vxe^3=Vxe$ and so $Ux=Uxe^2y\relrho Vxey=Vy$, which
  yields $U=Uxy\relrho Vy^2$. Thus $U\relrho Vy^2\relrho Uy^4$.

  From our observations about normal form words, $U$ is of the form
  $pe^\gamma q$, where $p$ is left-invertible, $q$ is
  right-invertible, and $\gamma\in\{0,1,2\}$. So we have
  $pe^{\gamma}q\relrho pe^{\gamma}qy^4$, which implies
  $e^{\gamma}q\relrho e^{\gamma}qy^4$ and thus, by left-multiplying by
  $e^{3-\gamma}$ and using $e^3=e$, we obtain $eq\relrho eqy^4$. Hence
  $eq\relrho eqy^{4k}$ for all $k\geq 1$. Since $q$ is
  right-invertible, it follows from the rewriting system and the form
  of right-invertible normal form words that for some appropriate
  $k\geq 1$ and $n \geq 1$ that $qy^{4k}=y^n$ and so $eq\relrho
  ey^n$. Then $y^n=xey^n\relrho xeq$. But $xeq$ is right-invertible in
  $M$ and so $y$ is right-invertible in $M/\relrho$. Together with
  $xy=1$, this implies that $y$, and thus $x$, and thus $e$, are
  invertible in $M/\relrho$ and so by the Lemma, $\relrho=M\times M$.

\item $\alpha=0$ and $\beta=1$. Then $Ux\relrho Vxe$. By right
  multiplying by $e$, we reduce to the previous case.

\item $\alpha=0$ and $\beta=2$. Then $Ux\relrho Vxe^2$ and so
  $U=Uxe^2y^2\relrho Vxe^2\cdot e^2y^2=Vxe^ey^2=V$, whence $Ux\relrho
  Uxe^2$. Then $U=Uxy\relrho Uxe^2y=Ux$ and so $U=Uxy\relrho Uy$. Now
  proceeding in a similar way to sub-case~a) (where we had $U
  \relrho Uy^4$) we obtain $\relrho=M\times M$.
\end{enumerate}

\item Only one of $u$ and $v$ contains letters $x$. Interchanging $u$
  and $v$ if necessary, assume that $u$ contains a letter $x$ and $v$
  does not; that is, $v\in\{e,y\}^*$. We have three cases to consider:
\begin{enumerate}
\item $u\equiv Ux$. Then $v\relrho Ux=Uxe^2y\relrho ve^2y$. Then
  $e\relrho e^3y=ey$ and so $y=xey\relrho xe$, whence $y^2\relrho
  xey=y$. Then $xy=1$ implies $y\relrho 1$ and so, proceeding as in
  the base case of the induction, we obtain $\relrho=M\times M$.

\item $u\equiv Uxe$. Then $v\relrho Uxe=Uxe^2ye\relrho veye$. Since
  $v$ is in normal form and also in $\{e,y\}^*$, we have $v \equiv
  v'e^\alpha$ where $v'$ is left-invertible and $\alpha \in
  \{0,1,2\}$. So $e^{\alpha}\relrho
  e^{\alpha+1}ye$. Left-multiplication by $e^{4-\alpha}$ gives Then $e
  = e^4 \relrho e^5ye = e^2ye$ and so $ye = xeye \relrho xe^2$. Then
  $yey^2\relrho xe^2y^2=1$, whence $y$ is right-invertible in
  $M/\relrho$. By $xy=1$, this implies that $y$, $x$, and $e$, are all
  invertible in $M/\relrho$ and so by the Lemma, $\relrho=M\times M$.

\item $u\equiv Uxe^2$. Then $v\relrho Uxe^2=Uxe^2ye^2\relrho vye^2$
  (using $x = xe^2y$). As in sub-case~b), we have $v \equiv
  v'e^\alpha$ where $v'$ is left-invertible and $\alpha \in
  \{0,1,2\}$. Thus $e^\alpha \relrho e^\alpha y
  e^2$. Left-multiplication by $e^{3-\alpha}$ gives $e = e^3 \relrho
  e^3ye^2 = eye^2$. Then $xe\relrho xeye^2=ye^2$ and so $x^2e\relrho
  e^2$, whence $1=xy = x^2ey\relrho e^2y$. Left-multiplication by $xe$
  gives $xe \relrho xe^3y = xey = y$, and so $ye^2 \relrho xe \relrho
  y$, and left-multiplication by $x$ gives $e^2 \relrho 1$. Hence $x =
  xe^2y \relrho xy = 1$ and so, proceeding as in the base case of the
  induction, we obtain $\relrho=M\times M$.
\end{enumerate}

\item Neither $u$ nor $v$ contain letters from $x$ but both contain
  letters from $y$. That is, $u\equiv e^{\alpha}yU$ and $v\equiv
  e^{\beta}yV$ for some $\alpha,\beta\in\{0,1,2\}$. If $\alpha=\beta$,
  then $U \neq V$ (since $u \neq v$) and $U=x^2e^{3-\alpha}u\relrho
  x^2e^{3-\alpha}v=V$ and so by induction $\rho = M \times M$ and the
  proof is complete. So, we may assume that $\alpha<\beta$. We have
  three natural cases:
\begin{enumerate}
\item $\alpha=1$ and $\beta=2$. Then $eyU\relrho e^2yV$ and so
  $yU=xeyU\relrho xe^2yV=xV$, whence $U=xyU\relrho x^2V$. Also
  $e^2yU\relrho eyV$ and therefore $xU=xe^2yU\relrho xeyV=yV$, whence
  $x^2U\relrho xyV=V$. Thus $U\relrho x^2V\relrho x^4U$. By reasoning
  symmetrical to sub-case~1(a), this leads to $\relrho=M\times M$.

\item $\alpha=0$ and $\beta=1$. By left multiplying by $e$, we reduce
  to the previous case.

\item $\alpha=0$ and $\beta=2$. Then $yV\relrho e^2yV$ and so $U = xyV
  \relrho xe^2yV = xV$. Also $eyU\relrho e^3yV=eyV$ and so
  $U=x^2eyU\relrho x^2eyV=V$, whence $U\relrho xU$ and by reasoning
  symmetrical to sub-case~1(c) we obtain $\relrho=M\times M$.
\end{enumerate}

\item Neither of $u$ and $v$ contain letters $x$ and only one contains
  a letter $y$. Interchanging $u$ and $v$ if necessary, assume that
  $u$ contains a letter $y$: that is, $u \equiv e^{\alpha}yU$. Then
  $eyU=e^{3-\alpha}u\relrho e^{3-\alpha}v$ and so $yU=xeyU\relrho
  xe^{3-\alpha}v$. But $xe^{3-\alpha}v$ is right-invertible in $M$,
  and so $y$ is right-invertible in $M/\relrho$ and so
  $\relrho=M\times M$.

\item Neither $u$ nor $v$ contains letters $x$ or $y$. That is, $u \equiv
  e^\alpha$ and $v \equiv e^\beta$ for $\alpha,\beta \in \{0,1,2\}$ with
  $\alpha \neq \beta$. We have three sub-cases:
\begin{enumerate}
 \item $\alpha=1$ and $\beta=2$. Then $x = xey \relrho xe^2y =
   y$. Together with $xy=1$, this implies that $x$ and $y$, and thus
   $e$ are invertible in $M/\relrho$. So by the Lemma, $\rho = M
   \times M$.
\item $\alpha = 0$ and $\beta=1$. By multiplying by $e$, we reduce to
  the previous case.
\item $\alpha = 0$ and $\beta=2$. Then $x = xe^2y \relrho xy = 1$ and
  proceeding as in the base case of the induction shows that $\rho = M
  \times M$.
\end{enumerate}
\end{enumerate}
This completes the induction step, and so $M$ is congruence-free.
\end{proof}


\end{document}